\documentclass{amsart}
\newtheorem{thm}{Theorem}
\newtheorem{lem}[thm]{Lemma}
\newtheorem{cor}[thm]{Corollary}

\theoremstyle{definition}

\newtheorem{ex}[thm]{Example}
\newtheorem{rem}[thm]{Remark}

\providecommand{\abs}[1]{\lvert#1\rvert}
\providecommand{\Abs}[1]{\Bigl\lvert#1\Bigr\rvert}
\providecommand{\norm}[1]{\lVert#1\rVert}
\date{11 August, 2022}

\newcommand\dto{\overset{dist}\longrightarrow}

\begin{document}
\title[CLT for $m$-dependent random variables]
{Quantitative bounds in the central limit theorem for $m$-dependent random variables}

\author{Svante Janson}
\address{Svante Janson, Department of Mathematics, Uppsala University, PO Box 480, SE-751~06 Uppsala, Sweden, URL: http://www.math.uu.se/svante-janson}
\email{svante.janson@math.uu.se}

\author{Luca Pratelli}
\address{Luca Pratelli, Accademia Navale, viale Italia 72, 57100 Livorno, Italy}
\email{pratel@mail.dm.unipi.it}

\author{Pietro Rigo}
\address{Pietro Rigo, Dipartimento di Scienze Statistiche ``P. Fortunati'', Universit\`a di Bologna, via delle Belle Arti 41, 40126 Bologna, Italy}
\email{pietro.rigo@unibo.it}

\keywords{Central limit theorem; Lindeberg condition; $m$-dependent random variables; Quantitative bound; Total variation distance; Wasserstein distance}

\subjclass[2020]{60F05}

\begin{abstract}
For each $n\ge 1$, let $X_{n,1},\ldots,X_{n,N_n}$ be real random variables and $S_n=\sum_{i=1}^{N_n}X_{n,i}$. Let $m_n\ge 1$ be an integer. Suppose $(X_{n,1},\ldots,X_{n,N_n})$ is $m_n$-dependent, $E(X_{ni})=0$, $E(X_{ni}^2)<\infty$ and $\sigma_n^2:=E(S_n^2)>0$ for all $n$ and $i$. Then,
\begin{gather*}
d_W\Bigl(\frac{S_n}{\sigma_n},\,Z\Bigr)\le 30\,\bigl\{c^{1/3}+12\,U_n(c/2)^{1/2}\bigr\}\quad\quad\text{for all }n\ge 1\text{ and }c>0,
\end{gather*}
where $d_W$ is Wasserstein distance, $Z$ a standard normal random variable and
$$U_n(c)=\frac{m_n}{\sigma_n^2}\,\sum_{i=1}^{N_n}E\Bigl[X_{n,i}^2\,1\bigl\{\abs{X_{n,i}}>c\,\sigma_n/m_n\bigr\}\Bigr].$$
Among other things, this estimate of $d_W\bigl(S_n/\sigma_n,\,Z\bigr)$ yields a similar estimate of $d_{TV}\bigl(S_n/\sigma_n,\,Z\bigr)$ where $d_{TV}$ is total variation distance.

\end{abstract}

\maketitle

\section{Introduction}\label{intro}

Central limit theorems (CLTs) for $m$-dependent random variables have a
long history. Pioneering results, for a fixed $m$, were given by Hoeffding and Robbins
\cite{HOEROB} and Diananda \cite{DIAN} (for $m$-dependent sequences), and Orey \cite{OREY} (more generally,
and also for triangular arrays). These results were then extended to the case of increasing
$m = m_n$; see e.g.\ Bergstr\"om \cite{BERG}, Berk \cite{BERK}, Rio
\cite{RIO}, Romano and Wolf \cite{RW}, and Utev \cite{UT1}, \cite{UT2}.

\medskip

Obviously, CLTs for $m$-dependent random variables are often corollaries of more general results obtained under mixing conditions. A number of CLTs under mixing conditions are actually available. Without any claim of being exhaustive, we mention \cite{BRAD}, \cite{DMR}, \cite{PEL}, \cite{RIO}, \cite{UT1}, \cite{UT2} and references therein. However, mixing conditions are not directly related to our purposes (as stated below) and they will not be discussed further.

\medskip

This paper deals with an $(m_n)$-dependent array of random variables, where $(m_n)$ is any sequence of integers, and provides an upper bound for the Wasserstein distance between the standard normal law and the distribution of the normalized partial sums. A related bound for the total variation distance is obtained as well. To be more precise, we need some notation.

\medskip

For each $n\ge 1$, let $1\le m_n\le N_n$ be integers, $(X_{n,1},\ldots,X_{n,N_n})$ a collection of real random variables, and
$$S_n=\sum_{i=1}^{N_n}X_{n,i}.$$
Suppose
\begin{gather}\label{c1}
(X_{n,1},\ldots,X_{n,N_n})\text{ is }m_n\text{-dependent for every }n,
\end{gather}
\begin{gather}\label{c2}
E(X_{ni})=0,\quad E(X_{ni}^2)<\infty,\quad\sigma_n^2:=E(S_n^2)>0\quad\text{for all }n\text{ and }i,
\end{gather}
and define
$$U_n(c)=\frac{m_n}{\sigma_n^2}\,\sum_{i=1}^{N_n}E\Bigl[X_{n,i}^2\,1\bigl\{\abs{X_{n,i}}>c\,\sigma_n/m_n\bigr\}\Bigr]\quad\quad\text{for all }c>0.$$
Our main result (Theorem \ref{m2}) is that
\begin{gather}\label{n7d5}
d_W\Bigl(\frac{S_n}{\sigma_n},\,Z\Bigr)\le 30\,\bigl\{c^{1/3}+12\,U_n(c/2)^{1/2}\bigr\}\quad\quad\text{for all }n\ge 1\text{ and }c>0,
\end{gather}
where $d_W$ is Wasserstein distance and $Z$ a standard normal random variable.

\medskip

Inequality \eqref{n7d5} provides a quantitative estimate of $d_W\bigl(S_n/\sigma_n,\,Z\bigr)$. The connections between \eqref{n7d5} and other analogous results are discussed in Remark \ref{hb88j9w} and Section \ref{tv21d}. To our knowledge, however, no similar estimate of $d_W\bigl(S_n/\sigma_n,\,Z\bigr)$ is available under conditions \eqref{c1}--\eqref{c2} only. In addition, inequality \eqref{n7d5} implies the following useful result:
\begin{thm}[Utev \cite{UT1,UT2}]\label{a9ny65vg}
$S_n/\sigma_n\overset{dist}\longrightarrow Z$ provided conditions \eqref{c1}--\eqref{c2} hold and $U_n(c)\rightarrow 0$ for every $c>0$.
\end{thm}

\medskip

Based on inequality \eqref{n7d5}, we also obtain quantitative bounds for $d_K\bigl(S_n/\sigma_n,\,Z\bigr)$ and $d_{TV}\bigl(S_n/\sigma_n,\,Z\bigr)$, where $d_K$ and $d_{TV}$ are Kolmogorov distance and total variation distance, respectively. As to $d_K$, it suffices to recall that
\begin{gather*}
d_K\Bigl(\frac{S_n}{\sigma_n},\,Z\Bigr)\le 2\,\sqrt{d_W\Bigl(\frac{S_n}{\sigma_n},\,Z\Bigr)};
\end{gather*}
see Lemma \ref{b66z3e0ju5}. To estimate $d_{TV}$, define
$$l_n=2\,\int_0^\infty t\,\abs{\phi_n(t)}\,dt$$
where $\phi_n$ is the characteristic function of $S_n/\sigma_n$. By a result
in \cite{PR} (see Theorem \ref{math187v4} below),
$$d_{TV}\Bigl(\frac{S_n}{\sigma_n},\,Z\Bigr)\le
2\,d_W\Bigl(\frac{S_n}{\sigma_n},\,Z\Bigr)^{1/2}+
l_n^{2/3}\,d_W\Bigl(\frac{S_n}{\sigma_n},\,Z\Bigr)^{1/3}.
$$
Hence, $d_{TV}\bigl(S_n/\sigma_n,\,Z\bigr)$ can be upper bounded via inequality \eqref{n7d5}. For instance, in addition to \eqref{c1}--\eqref{c2}, suppose $X_{ni}\in L_\infty$ for all $n$ and $i$ and define
$$c_n=\frac{2\,m_n}{\sigma_n}\,\max_i\norm{X_{ni}}_\infty.$$
On noting that $U_n(c_n/2)=0$, one obtains
\begin{gather*}
d_{TV}\Bigl(\frac{S_n}{\sigma_n},\,Z\Bigr)\le \sqrt{120}\,\,c_n^{1/6}+30^{1/3}\,l_n^{2/3}\,c_n^{1/9}.
\end{gather*}

\medskip

The rest of this paper is organized into three sections. Section
\ref{notdef} just recalls some definitions and known results, Section
\ref{was454c} is devoted to proving inequality \eqref{n7d5}, while Section
\ref{tv21d} investigates $d_{TV}\bigl(S_n/\sigma_n,\,Z\bigr)$ and the
convergence rate provided by \eqref{n7d5}.

\medskip

The numerical constants in our results are obviously not best possible;
we have not tried to optimize them.
More important are the powers, $c^{1/3}$ and $U_n(c/2)^{1/2}$ in  \eqref{n7d5}
and similar powers in other results; we do not believe that these are optimal.
This is discussed in Section \ref{tv21d}.
How far \eqref{n7d5} can be improved, however, is essentially an open problem.

\medskip

\section{Preliminaries}\label{notdef}

The same notation as in Section \ref{intro} is adopted in the sequel.
It is implicitly assumed  that, for each $n\ge 1$, the variables
$(X_{ni}:1\le i\le N_n)$ are defined
on the same probability space (which may depend on $n$).

\medskip

Let $k\ge 0$ be an integer. A (finite or infinite) sequence $(Y_i)$ of random variables
is $k$-dependent if $(Y_i:i\le j)$ is independent of $(Y_i:i>j+k)$ for every $j$. In particular, $0$-dependent is the same as
independent. Given a sequence $(k_n)$ of non-negative integers, an array $(Y_{ni}:n\ge 1,$ $1\le i\le N_n)$ is said to be $(k_n)$-dependent if $(Y_{ni}:1\le i\le N_n)$ is $k_n$-dependent for every $n$.

\medskip

Let $X$ and $Y$ be real random variables. Three well known distances between their probability distributions are Wasserstein's, Kologorov's and total variation. Kolmogorov distance and total variation distance are, respectively,
\begin{gather*}
d_K(X,Y)=\sup_{t\in\mathbb{R}}\,\abs{P(X\le t)-P(Y\le t)}\quad\text{and}
\\d_{TV}(X,Y)=\sup_{A\in\mathcal{B}(\mathbb{R})}\,\abs{P(X\in A)-P(Y\in A)}.
\end{gather*}
Under the assumption $E\abs{X}+E\abs{Y}<\infty$, Wasserstein distance is
$$d_W(X,Y)=\inf_{U\sim X, V\sim Y}\,E\abs{U-V}$$
where $\inf$ is over the real random variables $U$ and $V$, defined on the same probability space, such that $U\sim X$ and $V\sim Y$. Equivalently,
$$d_W(X,Y)=\int_{-\infty}^\infty\abs{P(X\le t)-P(Y\le t)}\,dt=\sup_f\,\abs{Ef(X)-Ef(Y)}$$
where $\sup$ is over the 1-Lipschitz functions
$f:\mathbb{R}\rightarrow\mathbb{R}$. The next lemma is certainly known, but
we give a proof since we do not know of any reference for the first claims.

\begin{lem}\label{b66z3e0ju5}
Suppose $EX^2\le 1$, $EY^2\le 1$ and $EY=0$. Then,
\begin{align*}
  d_W(X,Y)&\le\sqrt{2},
\\
  d_W(X,Y)&\le4 \sqrt{d_K(X,Y)}.
\end{align*}
If $Y\sim N(0,1)$, we also have
$$
d_K(X,Y)\le 2\,\sqrt{d_W(X,Y)}.$$
\end{lem}

\begin{proof}
Take $U$ independent of $V$ with $U\sim X$ and $V\sim Y$. Then,
$$d_W(X,Y)\le E\abs{U-V}\le
\bigl\{E\bigl((U-V)^2\bigr)\bigr\}^{1/2}\le\sqrt{2}.$$
Moreover, for each $c>0$,
\begin{align}\label{sj1a}
d_W(X,Y)
&\notag=\int_{-\infty}^\infty\abs{P(X\le t)-P(Y\le t)}\,dt
\\ \notag
&\notag \le 2\,c\,d_K(X,Y)+\int_c^\infty\abs{P(X>t)-P(Y>t)}\,dt
\\ \notag
& +\int_c^\infty\abs{P(-X>t)-P(-Y>t)}\,dt
\\ \notag
&\le 2\,c\,d_K(X,Y)+\int_c^{\infty}\bigl\{ P(\abs{X}>t)+P(\abs{Y}>t)\bigr\}\,dt
\\ \notag
&\le 2\,c\,d_K(X,Y)+\int_c^{\infty}\frac{2}{t^2}\,dt= 2\,c\,d_K(X,Y)+\frac{2}{c}.
 \end{align}
Hence, letting $c=d_K(X,Y)^{-1/2}$, one obtains $d_W(X,Y)\le 4\,\sqrt{d_K(X,Y)}$.

Finally, if $Y\sim N(0,1)$, it is well known that $d_K(X,Y)\le
2\,\sqrt{d_W(X,Y)}$; see e.g.\ \cite[Theorem 3.3]{CHEN}.
\end{proof}

\medskip

Finally, under some conditions, $d_{TV}$ can be estimated through $d_W$.
We report a result which allows this; in our setting we simply
take $V=1$ below.

\begin{thm}[A version of {\cite[Theorem 1]{PR}}]\label{math187v4}
Let $X_n,\,V,\,Z$ be real random variables,
and suppose that
$Z\sim N(0,1)$, $V>0$, $EV^2=EX_n^2=1$ for all $n$, and $V$ is
independent of $Z$.
Let
$\phi_n$ be the characteristic function of $X_n$, and
$$l_n=2\,\int_0^\infty t\,\abs{\phi_n(t)}\,dt.$$
Then,
$$d_{TV}(X_n,\,VZ)\le\bigl\{1+E(1/V)\bigr\}\,d_W(X_n,\,VZ)^{1/2}+
l_n^{2/3}\,d_W(X_n,\,VZ)^{1/3}$$
for each $n$.
\end{thm}

\begin{proof}
This is essentially a special case  of \cite[Theorem 1]{PR}, with
$\beta=2$ and the constant $k$ made explicit.
Also, the assumption $d_W(X_n,\,VZ)\rightarrow 0$ in \cite[Theorem 1]{PR}
is not needed; we use instead $d_W(X_n,\,VZ)\le\sqrt{2}$ from Lemma
\ref{b66z3e0ju5}.
Using this and $EX_n^2=1$,
the various constants
appearing in the proof can be explicitly evaluated. In fact, improving the argument in \cite{PR} slightly
by using $P(|X_n|>t)\le EX_n^2/t^2=t^{-2}$,
and as just said using $d_W(X_n,\,VZ)\le\sqrt{2}$,
we can take $k^*=5+4\sqrt2$ in the proof. After simple calculations, this implies that the constant
$k$ in \cite{PR} can be taken as
$$
k=\frac{1}{2}\cdot\frac{3}{2}\cdot 2^{1/3} (5+4\sqrt2)^{1/3}\pi^{-2/3}
<1.
$$
\end{proof}

\medskip

\section{An upper bound for Wasserstein distance}\label{was454c}

As noted in Section \ref{intro}, our main result is:

\begin{thm}\label{m2}
Under conditions \eqref{c1}--\eqref{c2},
$$d_W\Bigl(\frac{S_n}{\sigma_n},\,Z\Bigr)\le 30\,\bigl\{c^{1/3}+12\,U_n(c/2)^{1/2}\bigr\}$$
for all $n\ge 1$ and $c>0$, where $Z$ denotes a standard normal random variable.
\end{thm}

\medskip

In turn, Theorem \ref{m2} follows from the following result,
which is a sharper version of the special case $m_n=1$.

\medskip

\begin{thm}\label{m1}
Let $X_1,\ldots,X_N$ be real random variables and $S=\sum_{i=1}^NX_i$. Suppose $(X_1,\ldots,X_N)$ is $1$-dependent and
\begin{gather*}
E(X_i)=0,\,\,E(X_i^2)<\infty\text{ for all }i\text{ and }\sigma^2:=E(S^2)>0.
\end{gather*}
Then,
\begin{gather*}
d_W\Bigl(\frac{S}{\sigma},\,Z\Bigr)\le 30\,\bigl\{c^{1/3}+6\,L(c)^{1/2}\bigr\}\quad\quad\text{for all }c>0,
\end{gather*}
where $Z$ is a standard normal random variable and
$$L(c)=\frac{1}{\sigma^2}\,\sum_{i=1}^NE\Bigl[X_i^2\,1\bigl\{\abs{X_i}>c\,\sigma\bigr\}\Bigr].$$
\end{thm}

\medskip

To deduce Theorem \ref{m2} from Theorem \ref{m1}, define $M_n=\lceil N_n/m_n\rceil$, $X_{n,i}=0$ for $i>N_n$, and
\begin{gather*}
Y_{n,i}=\sum_{j=(i-1)m_n+1}^{im_n}X_{n,j}\quad\quad\text{for }i=1,\ldots,M_n.
\end{gather*}
Since $(Y_{n,1},\ldots,Y_{n,M_n})$ is 1-dependent and $\sum_iY_{n,i}=\sum_iX_{n,i}=S_n$, Theorem \ref{m1} implies
\begin{gather}\label{ut4r5bh}
d_W\Bigl(\frac{S_n}{\sigma_n},\,Z\Bigr)\le 30\,\bigl\{c^{1/3}+6\,L_n(c)^{1/2}\bigr\}
\end{gather}
where
\begin{gather*}
L_n(c)=\frac{1}{\sigma_n^2}\,\sum_{i=1}^{M_n}E\Bigl[Y_{n,i}^2\,1\bigl\{\abs{Y_{n,i}}>c\,\sigma_n\bigr\}\Bigr].
\end{gather*}
Therefore, to obtain Theorem \ref{m2}, it suffices to note the following
inequality:
\begin{lem}\label{Llu}
  With notations as above,
for every $c>0$,
\begin{gather*}
L_n(2c)\le 4\,U_n(c).
\end{gather*}
\end{lem}

\medskip

In the rest of this section, we prove Lemma \ref{Llu} and Theorem \ref{m1}.
We also obtain a (very small) improvement of Utev's Theorem \ref{a9ny65vg}.

\subsection{Proof of Lemma \ref{Llu} and Utev's theorem}
\begin{proof}[Proof of Lemma \ref{Llu}]
 Fix $c>0$ and define
\begin{gather*}
V_{n,i}=\sum_{j=(i-1)m_n+1}^{im_n}X_{n,j}\,1\{\abs{X_{n,j}}>c\,\sigma_n/m_n\}.
\end{gather*}
Since $\abs{Y_{n,i}}\le\abs{V_{n,i}}+c\,\sigma_n$, one obtains
\begin{gather*}
\abs{Y_{n,i}}\,1\{\abs{Y_{n,i}}>2\,c\,\sigma_n\}\le\bigl(\abs{V_{n,i}}+c\,\sigma_n\bigr)\,1\{\abs{V_{n,i}}>c\,\sigma_n\}\le 2\,\abs{V_{n,i}}.
\end{gather*}
Therefore,
\begin{align*}
\sigma_n^2\,L_n(2c) &=\sum_{i=1}^{M_n}E\bigl[Y_{n,i}^2\,1\{\abs{Y_{n,i}}>2\,c\,\sigma_n\}\bigr]\le 4\,\sum_{i=1}^{M_n}E(V_{n,i}^2)
\\
&\le 4\,m_n\,\sum_{i=1}^{M_n}\sum_{j=(i-1)m_n+1}^{im_n}E\bigl[X_{n,j}^2\,1\{\abs{X_{n,j}}>c\,\sigma_n/m_n\}\bigr]
\\
&=4\,m_n\,\sum_{i=1}^{N_n}E\bigl[X_{n,i}^2\,1\{\abs{X_{n,i}}>c\,\sigma_n/m_n\}\bigr]=4\sigma_n^2\,U_n(c).
\end{align*}
\end{proof}

\medskip

We also note that, because of \eqref{ut4r5bh}, Theorem \ref{m1} implies:
\begin{cor}\label{z4oao7d}
$S_n/\sigma_n\overset{dist}\longrightarrow Z$ if conditions \eqref{c1}--\eqref{c2} hold and $L_n(c)\rightarrow 0$ for every $c>0$.
\end{cor}
Corollary \ref{z4oao7d} slightly improves Theorem \ref{a9ny65vg}. In fact,
$U_n(c)\rightarrow 0$ for all $c>0$
implies $L_n(c)\rightarrow 0$ for all $c>0$, because of Lemma \ref{Llu}, but
the converse is not true.

\begin{ex}\textbf{($L_n(c)\rightarrow 0$ does not imply $U_n(c)\rightarrow 0$).}
Let $(V_n:n\ge 1)$ be an i.i.d.\ sequence of real random variables such that $V_1$ is absolutely continuous with density $f(x)=(3/2)\,x^{-4}\,1_{[1,\infty)}(\abs{x})$. Let $m_n$ and $t_n$ be positive integers such that $m_n\rightarrow\infty$. Define $N_n=m_n\,(t_n+1)$ and
\begin{gather*}
X_{n,i}=V_i\,\text{ if }1\le i\le m_nt_n\quad\text{and}\quad X_{n,i}=V_{m_nt_n+1}\,\text{ if }m_nt_n<i\le m_n(t_n+1).
\end{gather*}
Define also
\begin{gather*}
T_n=\frac{\sum_{j=1}^{m_n}V_j}{\sqrt{m_n}}.
\end{gather*}
Then, $EV_1^2=3$, $\sigma_n^2=3\,(m_nt_n+m_n^2)$ and
\begin{align*}
L_n(c)&=\frac{1}{\sigma_n^2}\,\sum_{i=1}^{M_n}E\bigl[Y_{n,i}^2\,1\{\abs{Y_{n,i}}>c\,\sigma_n\}\bigr]\le\frac{1}{\sigma_n^2}
\,\sum_{i=1}^{t_n}E\bigl[Y_{n,i}^2\,1\{\abs{Y_{n,i}}>c\,\sigma_n\}\bigr]+\frac{3\,m_n^2}{\sigma_n^2}
\\
&=\frac{m_nt_n}{\sigma_n^2}\,E\bigl[T_n^2\,1\{\abs{T_n}>c\,\sigma_n/\sqrt{m_n}\}\bigr]+\frac{3\,m_n^2}{\sigma_n^2}.
\end{align*}
If $m_n=\,$o$(t_n)$, then $m_n^2/\sigma_n^2\rightarrow 0$,
$m_nt_n/\sigma_n^2\rightarrow 1/3$ and
$\sigma_n/\sqrt{m_n}\rightarrow\infty$. Moreover, the sequence $(T_n^2)$ is
uniformly integrable
(since $T_n\dto N(0,3)$ with (trivial) convergence of second moments). Hence, if $m_n=\,$o$(t_n)$, one obtains, for every $c>0$,
\begin{gather*}
\limsup_n L_n(c)\le\frac{1}{3}\,\limsup_nE\bigl[T_n^2\,1\{\abs{T_n}>c\,\sigma_n/\sqrt{m_n}\}\bigr]=0.
\end{gather*}
However,
\begin{align*}
U_n(c)&=\frac{m_n}{\sigma_n^2}\,\sum_{i=1}^{N_n}E\bigl[X_{n,i}^2\,1\{\abs{X_{n,i}}>c\,\sigma_n/m_n\}\bigr]=\frac{m_nN_n}{\sigma_n^2}\,E\bigl[V_1^2\,1\{\abs{V_1}>c\,\sigma_n/m_n\}\bigr]
\\
&=\frac{3\,m_nN_n}{\sigma_n^2}\,\int_{c\,\sigma_n/m_n}^\infty x^{-2}dx
=\frac{3N_n}{c\,\sigma_n^2}\,\frac{m_n^2}{\sigma_n}
\ge \frac{3 t_n m_n^3}{c(6m_nt_n)^{3/2}}
\end{align*}
for each $n$ such that $c\,\sigma_n/m_n\ge 1$ and $m_n\le t_n$.
Therefore, $L_n(c)\rightarrow 0$ and $U_n(c)\rightarrow\infty$
for all $c>0$
whenever $m_n=\,$o$(t_n)$ and $t_n=\,$o$(m_n^3)$. This happens, for instance, if $m_n\rightarrow\infty$ and $t_n=m_n^2$.
\end{ex}

\medskip

\subsection{Proof of Theorem \ref{m1}}
Our proof of Theorem \ref{m1} requires three lemmas. A result by R\"ollin \cite{ROLLIN} plays a crucial role in one of them (Lemma \ref{vb1a0i}).

\medskip

In this subsection, $X_1,\ldots,X_N$ are real random variables and $S=\sum_{i=1}^NX_i$. We assume that $(X_1,\ldots,X_N)$ is 1-dependent and
\begin{gather*}
E(X_i)=0,\,\,E(X_i^2)<\infty\text{ for all }i\text{ and }\sigma^2:=E(S^2)>0.
\end{gather*}
Moreover, $Z$ is a standard normal random variable {\em independent of} $(X_1,\ldots,X_N)$.

\medskip

For each $i=1,\ldots,N$, define
\begin{gather*}
Y_i=X_i-E(X_i\mid\mathcal{F}_{i-1})+E(X_{i+1}\mid\mathcal{F}_i)
\end{gather*}
where $\mathcal{F}_0$ is the trivial $\sigma$-field, $\mathcal{F}_i=\sigma(X_1,\ldots,X_i)$ and $X_{N+1}=0$. Then,
$$
E(Y_i\mid\mathcal{F}_{i-1})=0\text{ for all }i\text{ and
}\sum_{i=1}^NY_i=\sum_{i=1}^NX_i=S
\text{ a.s.}
$$

\begin{lem}\label{vg8uj2}
Let $\gamma>0$ be a constant and
$V^2=\sum_{i=1}^NE(Y_i^2\mid\mathcal{F}_{i-1})$. Then,
$$E\Bigl\{\Bigl(\frac{V^2}{\sigma^2}-1\Bigr)^2\Bigr\}\le 16\,\gamma^2$$
provided $\max_i\abs{X_i}\le\sigma\,\gamma/3$ a.s.
\end{lem}

\begin{proof}
First note that
$$\sigma^2=E(S^2)=E\bigl\{\bigl(\sum_{i=1}^NY_i\bigr)^2\bigr\}=\sum_{i=1}^NE(Y_i^2)=E\bigl(\sum_{i=1}^NY_i^2\bigr).$$
Moreover, since $\max_i\abs{Y_i}\le\gamma\,\sigma$ a.s., one obtains
\begin{gather*}
\sum_{i=1}^NE(Y_i^4)\le\gamma^2\sigma^2\,\sum_{i=1}^NE(Y_i^2)=\gamma^2\sigma^4.
\end{gather*}
Therefore,
\begin{align*}
E\Bigl\{\Bigl(\frac{V^2}{\sigma^2}-1\Bigr)^2\Bigr\}&\le\frac{2}{\sigma^4}\,\Bigl\{E\Bigl[\Bigl(\,\sum_{i=1}^N(E(Y_i^2\mid\mathcal{F}_{i-1})-Y_i^2)\Bigr)^2\Bigr]+\text{Var}\Bigl(\,\sum_{i=1}^NY_i^2\Bigr)\Bigr\}
\\
&=\frac{2}{\sigma^4}\,\Bigl\{\sum_{i=1}^NE\Bigl(Y_i^4-E(Y_i^2\mid\mathcal{F}_{i-1})^2\Bigr)+\sum_{i=1}^N\text{Var}(Y_i^2)
\\ &\qquad +2\sum_{1\le i<j\le N}\text{Cov}(Y_i^2,\,Y_j^2)\Bigr\}\\ &\le\frac{4}{\sigma^4}\,\Bigl\{\sum_{i=1}^NE(Y_i^4)+\sum_{1\le i<j\le N}\text{Cov}(Y_i^2,\,Y_j^2)\Bigr\}\\
&\le 4\,\gamma^2+\frac{4}{\sigma^4}\,\sum_{1\le i<j\le N}\text{Cov}(Y_i^2,\,Y_j^2).
\end{align*}
To estimate the covariance part, define
$$Q_i=Y_i^2-E(Y_i^2)\quad\text{and}\quad T_i=\sum_{k=1}^iY_k=\sum_{k=1}^iX_k+E(X_{i+1}\mid\mathcal{F}_i).$$
For each fixed $1\le i<N$,
since $(T_1,\ldots,T_N)$ is a martingale,
\begin{align*}
\sum_{j>i}\text{Cov}(Y_i^2,\,Y_j^2)&=\sum_{j>i}E\bigl(Q_iY_j^2\bigr)=E\Bigl\{Q_i\,\sum_{j>i}Y_j^2\Bigr\}=E\Bigl\{Q_i\,\bigl(T_N-T_i\bigr)^2\Bigr\}
\\
&=E\Bigl\{Q_i\,\bigl(T_N-T_{i+1}\bigr)^2\Bigr\}
+E\bigl(Q_iY_{i+1}^2\bigr)
\\ &\le E\Bigl\{Q_i\,\bigl(T_N-T_{i+1}\bigr)^2\Bigr\}+E(Y_i^4)+E(Y_{i+1}^4).
\end{align*}
Finally, since $(X_1,\ldots,X_N)$ is 1-dependent,
$E Q_i=0$ and $E X_j=0$,
\begin{align*}
E\Bigl\{Q_i\,\bigl(T_N-T_{i+1}\bigr)^2\Bigr\}&=E\Bigl\{Q_i\,\Bigl(\sum_{k=i+2}^NX_k-E(X_{i+2}\mid\mathcal{F}_{i+1})\Bigr)^2\Bigr\}
\\
&=E\Bigl\{Q_i\,\Bigl(E(X_{i+2}\mid\mathcal{F}_{i+1})^2-2\,X_{i+2}\,E(X_{i+2}\mid\mathcal{F}_{i+1})\Bigr)\Bigr\}
\\
&=-E\Bigl\{Q_i\,E(X_{i+2}\mid
\mathcal{F}_{i+1})^2\Bigr\}\\
&\le E(Y_i^2)\,E\bigl\{E(X_{i+2}\mid\mathcal{F}_{i+1})^2\bigr\}\le\gamma^2\sigma^2E(Y_i^2).
\end{align*}
To sum up,
\begin{gather*}
E\Bigl\{\Bigl(\frac{V^2}{\sigma^2}-1\Bigr)^2\Bigr\}\le 4\,\gamma^2+\frac{4}{\sigma^4}\,\sum_{i=1}^{N-1}\Bigl(E(Y_i^4)+E(Y_{i+1}^4)+\gamma^2\sigma^2E(Y_i^2)\Bigr)\le 16\,\gamma^2.
\end{gather*}
\end{proof}

\begin{lem}\label{vb1a0i}
If $\max_i\abs{X_i}\le\sigma\,\gamma/3$ a.s., then
$$
d_W\Bigl(\frac{S}{\sigma}\,,\,Z\Bigr)
\le 16\,\gamma^{1/3}.
$$
\end{lem}

\begin{proof}
By Lemma \ref{b66z3e0ju5}, $d_W\bigl(S/\sigma,\,Z\bigr)\le\sqrt{2}$. Hence,
it can be assumed that $\gamma\le 1$.

Define
\begin{gather*}\tau=\max\bigl\{m:1\le m\le N,\,\sum_{k=1}^mE(Y_k^2/\sigma^2\mid\mathcal{F}_{k-1})\le 1\bigr\},
\\J_i=1\{\tau\ge i\}\,\frac{Y_i}{\sigma}+1\{\tau=i-1\}\,\Bigl(1-\sum_{k=1}^{i-1} E(Y_k^2/\sigma^2\mid\mathcal{F}_{k-1})\Bigr)^{1/2}Z\quad\text{for }i=1,\ldots,N,
\\J_{N+1}=1\{\tau=N\}\,\Bigl(1-\sum_{k=1}^N E(Y_k^2/\sigma^2\mid\mathcal{F}_{k-1})\Bigr)^{1/2}Z.
\end{gather*}
Since $\tau$ is a stopping time,
$Z$ is independent of $(X_1,\ldots,X_N)$,
and $E(Y_i\mid\mathcal{F}_{i-1})=0$,
one obtains
$$E(J_i\mid\mathcal{F}_{i-1})=0\text{ for all }i\text{ and }\sum_{k=1}^{N+1}E(J_k^2\mid\mathcal{F}_{k-1})=1\text{ a.s.}$$
Therefore, for each $a>0$, a result by R\"ollin \cite[Theorem 2.1]{ROLLIN} implies
$$d_W\Bigl(\sum_{i=1}^{N+1}J_i\,,\,Z\Bigr)\le 2a+\frac{3}{a^2}\,\sum_{i=1}^{N+1}E\abs{J_i}^3.$$
To estimate $E\abs{J_i}^3$ for $i\le N$, note that $E\abs{Z}^3\le 2$ and
$(1/\sigma)\,\max_i\abs{Y_i}\le\gamma$ a.s. Therefore,
for $1\le i\le N$,
\begin{align*}
E\abs{J_i}^3 &=E\Bigl\{1\{\tau\ge i\}\,\frac{\abs{Y_i}^3}{\sigma^3}\Bigr\}\kern-0.6mm+\kern-0.6mmE\Bigl\{1\{\tau=i-1\}\,\Bigl(1\kern-0.4mm-\kern-0.4mm\sum_{k=1}^{i-1} E(Y_k^2/\sigma^2\mid\mathcal{F}_{k-1})\Bigr)^{3/2}\abs{Z}^3\Bigr\}
\\
&\le\gamma\,E\Bigl\{1\{\tau\ge i\}\,\frac{Y_i^2}{\sigma^2}\Bigr\}
\\ &\qquad
+E\Bigl\{1\{\tau=i-1\}\,\Bigl(1-\sum_{k=1}^{i-1} E(Y_k^2/\sigma^2\mid\mathcal{F}_{k-1})\Bigr)^{1/2}\Bigr\}\,E\abs{Z}^3
\\
&\le\gamma\,E\Bigl\{1\{\tau\ge i\}\,\frac{Y_i^2}{\sigma^2}\Bigr\}
\\ &\qquad+2\,E\Bigl\{1\{\tau=i-1\}\,\Bigl(\sum_{k=1}^i E(Y_k^2/\sigma^2\mid\mathcal{F}_{k-1})-\sum_{k=1}^{i-1} E(Y_k^2/\sigma^2\mid\mathcal{F}_{k-1})\Bigr)^{1/2}\Bigr\}
\\
&=\gamma\,E\Bigl\{1\{\tau\ge i\}\,\frac{Y_i^2}{\sigma^2}\Bigr\}+2\,E\Bigl\{1\{\tau=i-1\}\,E(Y_i^2/\sigma^2\mid\mathcal{F}_{i-1})^{1/2}\Bigr\}
\\& \le\gamma\,E\Bigl\{1\{\tau\ge i\}\,\frac{Y_i^2}{\sigma^2}\Bigr\}+2\,\gamma\,P(\tau=i-1).
\end{align*}
Hence,
\begin{gather*}
\sum_{i=1}^NE\abs{J_i}^3
\le\gamma\,E\Bigl[\sum_{i=1}^N\frac{Y_i^2}{\sigma^2}\Bigr]+2\gamma
= 3\gamma.
\end{gather*}
Similarly,
\begin{align*}
E\abs{J_{N+1}}^3&=E\Bigl\{1\{\tau=N\}\,\Bigl(1-\sum_{k=1}^N E(Y_k^2/\sigma^2\mid\mathcal{F}_{k-1})\Bigr)^{3/2}\Bigr\}\,E\abs{Z}^3
\\
&\le 2\,E\Bigl\{1\{\tau=N\}\,\Bigl(1-\sum_{k=1}^N E(Y_k^2/\sigma^2\mid\mathcal{F}_{k-1})\Bigr)\Bigr\}
\\&\le 2\,E\Bigl\{\Bigl(1-\sum_{k=1}^N E(Y_k^2/\sigma^2\mid
\mathcal{F}_{k-1})\Bigr)^2\Bigr\}^{1/2}
\\
&=2\,E\Bigl\{\Bigl(1-\frac{V^2}{\sigma^2}\Bigr)^2\Bigr\}^{1/2}\le 8\,\gamma
\end{align*}
where the last inequality is due to Lemma \ref{vg8uj2}. It follows that
\begin{gather*}
d_W\Bigl(\sum_{i=1}^{N+1}J_i\,,\,Z\Bigr)
\le 2a+\frac{3}{a^2}(3\gamma + 8\gamma)=2a+\frac{33\,\gamma}{a^2}.
\end{gather*}

Next, we estimate $d_W\bigl(S/\sigma,\,\sum_{i=1}^{N}J_i\bigr)$. To this end, we let
$$W_i=\sum_{k=1}^i E(Y_k^2/\sigma^2\mid\mathcal{F}_{k-1})$$
and we note that
\begin{align*}
\frac{S}{\sigma}-\sum_{i=1}^{N}J_i&=\sum_{i=1}^{N}\Bigl(\frac{Y_i}{\sigma}-J_i\Bigr)=\sum_{i=1}^{N}1\{\tau<i\}\,\Bigl(\frac{Y_i}{\sigma}-J_i\Bigr)
\\&=\sum_{i=1}^{N-1}1\{\tau=i\}\,\Bigl\{\,\sum_{k=i+1}^N\frac{Y_k}{\sigma}-\bigl(1-W_i\bigr)^{1/2}Z\Bigr\}.
\end{align*}
Therefore, recalling the definition of $\tau$,
\begin{align*}
d_W\Bigl(\frac{S}{\sigma}\,,\,\sum_{i=1}^{N}J_i\Bigr)^2&\le \Bigl(E\Abs{\frac{S}{\sigma}-\sum_{i=1}^{N}J_i}\Bigr)^2\le E\Bigl\{\Bigl(\frac{S}{\sigma}-\sum_{i=1}^{N}J_i\Bigr)^2\Bigl\}
\\
&=\sum_{i=1}^{N-1}E\Bigl\{1\{\tau=i\}\,\Bigl\{\,\sum_{k=i+1}^N\frac{Y_k}{\sigma}-\bigl(1-W_i\bigr)^{1/2}Z\Bigr\}^2\Bigr\}
\\
&=\sum_{i=1}^{N-1}E\Bigl\{1\{\tau=i\}\,\Bigl\{\,\sum_{k=i+1}^N
E(Y_k^2/\sigma^2\mid
\mathcal{F}_{k-1})+1-W_i\Bigr\}\Bigr\}
\\
&\le\sum_{i=1}^{N-1}E\Bigl\{1\{\tau=i\}\,\Bigl\{\,\sum_{k=i+2}^NE(Y_k^2/\sigma^2\mid\mathcal{F}_{k-1})+2\,E(Y_{i+1}^2/\sigma^2\mid\mathcal{F}_i)\Bigr\}\Bigr\}
\\ &\le\sum_{i=1}^{N-1}E\Bigl\{1\{\tau=i\}\,\Bigl\{V^2/\sigma^2-1+2\,\gamma^2\Bigr\}\Bigr\}
\le E\abs{V^2/\sigma^2-1}+2\,\gamma^2\\&\le 4\,\gamma+2\,\gamma^2
\end{align*}
where the last inequality is because of Lemma \ref{vg8uj2}.
Since we assumed $\gamma\le 1$, we obtain
\begin{gather*}
d_W\Bigl(\frac{S}{\sigma}\,,\,\sum_{i=1}^{N}J_i\Bigr)\le \sqrt{6\,\gamma}.
\end{gather*}

Finally, using Lemma \ref{vg8uj2} again, one obtains
$$
d_W\Bigl(\sum_{i=1}^{N}J_i\,,\,\sum_{i=1}^{N+1}J_i\Bigr)\le
E\abs{J_{N+1}}\le E\Bigl\{\Abs{\frac{V^2}{\sigma^2}-1}^{1/2}\Bigr\}\le E\Bigl\{\Bigl(\frac{V^2}{\sigma^2}-1\Bigr)^2\Bigr\}^{1/4}\le 2\,\sqrt{\gamma}.$$
Collecting all these facts together,
\begin{align*}
d_W\Bigl(\frac{S}{\sigma}\,,\,Z\Bigr)&\le d_W\Bigl(\frac{S}{\sigma}\,,\,\sum_{i=1}^{N}J_i\Bigr)+d_W\Bigl(\sum_{i=1}^{N}J_i\,,\,\sum_{i=1}^{N+1}J_i\Bigr)+d_W\Bigl(\sum_{i=1}^{N+1}J_i\,,\,Z\Bigr)
\\
&\le \sqrt{6\,\gamma}+2\sqrt\gamma+2a+\frac{33\,\gamma}{a^2}\le 5\,\sqrt{\gamma}+2a+\frac{33\,\gamma}{a^2}.
\end{align*}
For $a=4\gamma^{1/3}$, the above inequality yields, using again $\gamma\le 1$,
\begin{gather*}
d_W\Bigl(\frac{S}{\sigma}\,,\,Z\Bigr)
\le 5\,\sqrt{\gamma}+\Bigl(8+\frac{33}{16}\Bigr)\,\gamma^{1/3}
\le 16 \,\gamma^{1/3}.
\end{gather*}
This concludes the proof.
\end{proof}

\medskip

\begin{rem}\label{hb88j9w}
If we do not care about the value of the constant in the estimate,
the proof of Lemma \ref{vb1a0i} could be shortened by exploiting a result by
Fan and Ma \cite{FANMA}; this result, however, does not provide
explicit values of the majorizing constants.
We also note that, under the conditions of Lemma \ref{vb1a0i}, Heyde--Brown's inequality \cite{HB} yields
$$d_K\Bigl(\frac{S}{\sigma}\,,\,Z\Bigr)\le b\,\Bigl\{E\Bigl(\Bigl(\frac{V^2}{\sigma^2}-1\Bigr)^2\Bigr)+\frac{1}{\sigma^4}\,\sum_{i=1}^NEY_i^4\Bigr\}^{1/5}$$
for some constant $b$ independent of $N$. By Lemmas \ref{b66z3e0ju5} and \ref{vg8uj2}, this implies
$$d_W\Bigl(\frac{S}{\sigma}\,,\,Z\Bigr)\le 4\,\sqrt{d_K\Bigl(\frac{S}{\sigma}\,,\,Z\Bigr)}\le 4\,\sqrt{b}\,\Bigl\{16\,\gamma^2+\,\frac{\gamma^2}{\sigma^2}\,\sum_{i=1}^NEY_i^2\Bigr\}^{1/10}=4\,\sqrt{b}\,17^{1/10}\,\gamma^{1/5}.$$
Hence, in this case, Lemma \ref{vb1a0i} works better than Heyde--Brown's inequality to estimate $d_W\bigl(S/\sigma\,,\,Z\bigr)$.
\end{rem}

\medskip
Recall $L(c)$ defined in Theorem \ref{m1}.
\begin{lem}\label{bhu8y}
Letting
$\sigma_c^2=\operatorname{Var}\Bigl(\sum_{i=1}^N\frac{X_i}{\sigma}\,1\bigl\{\abs{X_i}\le
c\sigma\bigr\}\Bigr)$, we have
$$
\abs{\sigma_c-1}
\le \abs{\sigma_c^2-1}
\le 13\,L(c)\quad\quad\text{for all }c>0.$$
\end{lem}

\begin{proof}
Fix $c>0$ and define
\begin{gather*}
A_i=\bigl\{\abs{X_i}> c\sigma\bigr\},\quad T_i=\frac{X_i}{\sigma}\,1_{A_i}-E\Bigl(\frac{X_i}{\sigma}\,1_{A_i}\Bigr),\quad V_i=\frac{X_i}{\sigma}\,1_{A_i^c}-E\Bigl(\frac{X_i}{\sigma}\,1_{A_i^c}\Bigr).
\end{gather*}
On noting that $\sigma_c^2=\text{Var}\,\Bigl(\,\sum_{i=1}^NV_i\Bigr)$, one obtains
$$1=\text{Var}\,\Bigl(\,\sum_{i=1}^N(T_i+V_i)\Bigr)=\text{Var}\,\Bigl(\,\sum_{i=1}^NT_i\Bigr)+\sigma_c^2+2\,\text{Cov}\,\Bigl(\sum_{i=1}^NT_i\,,\,\sum_{i=1}^NV_i\Bigr).$$
Since $(X_1,\ldots,X_N)$ is 1-dependent, it follows that
\begin{align*}
\abs{\sigma_c^2 -1}&\le\text{Var}\,\Bigl(\,\sum_{i=1}^NT_i\Bigr)+2\,\Abs{\text{Cov}\,\Bigl(\sum_{i=1}^NT_i\,,\,\sum_{i=1}^NV_i\Bigr)}
\\
&=\text{Var}\,\Bigl(\,\sum_{i=1}^NT_i\Bigr)+2\,\Abs{\sum_{i=1}^N\text{Cov}\,(T_i,V_i)
\\ &\qquad+\sum_{i=1}^{N-1}\text{Cov}\,(T_i,V_{i+1})+\sum_{i=2}^N\text{Cov}\,(T_i,V_{i-1})}.
\end{align*}
Moreover,
\begin{align}\label{sjb}
\text{Var}\,\Bigl(\,\sum_{i=1}^NT_i\Bigr)
&=\sum_{i=1}^N\text{Var}(T_i)+2\,\sum_{i=1}^{N-1}\text{Cov}(T_i,T_{i+1})
\\
&\notag
\le
\sum_{i=1}^N\text{Var}(T_i)
+\sum_{i=1}^{N-1}\Bigl(\text{Var}(T_i)+\text{Var}(T_{i+1})\Bigr)
\le 3\,L(c).
\end{align}
Similarly,
\begin{gather*}
\text{Cov}\,(T_i,V_i)=-E\Bigl(\frac{X_i}{\sigma}\,1_{A_i}\Bigr)\,E\Bigl(\frac{X_i}{\sigma}\,1_{A_i^c}\Bigr)=E\Bigl(\frac{X_i}{\sigma}\,1_{A_i}\Bigr)^2\le E\Bigl(\frac{X_i^2}{\sigma^2}\,1_{A_i}\Bigr)
\end{gather*}
and
\begin{align*}
\Abs{\text{Cov}\,(T_i,V_{i-1})}&\leq E\Bigl(\frac{\abs{X_iX_{i-1}}}{\sigma^2}\,1_{A_i}\,1_{A_{i-1}^c}\Bigr)+E\Bigl(\frac{\abs{X_i}}{\sigma}\,1_{A_i}\Bigr)\,E\Bigl(\frac{\abs{X_{i-1}}}{\sigma}\,1_{A_{i-1}^c}\Bigr)
\\
&\le 2\,c\,E\Bigl(\frac{\abs{X_i}}{\sigma}\,1_{A_i}\Bigr)\le 2\,E\Bigl(\frac{X_i^2}{\sigma^2}\,1_{A_i}\Bigr)
\end{align*}
where the last inequality is because
$$\frac{c\,\abs{X_i}}{\sigma}\,1_{A_i}\le\frac{\abs{X_i^2}}{\sigma^2}\,1_{A_i}.$$
By the same argument, $\Abs{\text{Cov}\,(T_i,V_{i+1})}\le 2\,\sigma^{-2}E\bigl(X_i^2\,1_{A_i}\bigr)$. Collecting all these facts together, one finally obtains
\begin{gather*}
\abs{\sigma_c^2 -1}\le 3\,L(c)+10\,\sum_{i=1}^NE\Bigl(\frac{X_i^2}{\sigma^2}\,1_{A_i}\Bigr)=13\,L(c).
\end{gather*}
This completes the proof, since obviously
$\abs{\sigma_c -1}\le\abs{\sigma_c^2 -1}$.
\end{proof}

\medskip

Having proved the previous lemmas, we are now ready to attack Theorem \ref{m1}.

\medskip
\begin{proof}[Proof of Theorem \ref{m1}]
Fix $c>0$. We have to show that
$$d_W\Bigl(\frac{S}{\sigma},\,Z\Bigr)\le 30\,\bigl\{c^{1/3}+6\,L(c)^{1/2}\bigr\}.$$
Since $d_W\bigl(S/\sigma,Z\bigr)\le\sqrt{2}$,
this inequality is trivially true if $L(c)\ge 1/100$ or if $c\ge 1$. Hence,
it can be assumed $L(c)<1/100$ and $c< 1$.
Then, Lemma \ref{bhu8y} implies $\sigma_c>0$.

Define $T_i$ and $V_i$ as in the proof of Lemma \ref{bhu8y}.
Then $|V_i|\le 2c$ for every $i$, and thus
$(V_1,\ldots,V_N)$
satisfies the conditions of Lemma \ref{vb1a0i} with
$\sigma$ replaced by $\sigma_c$ and
$\gamma=6\,c/\sigma_c$. Hence,
$$d_W\Bigl(\frac{\sum_{i=1}^NV_i}{\sigma_c}\,,\,Z\Bigr)
\le 16\,\bigl(6c/\sigma_c\bigr)^{1/3}.$$
Now, recall from \eqref{sjb}
that Var$\bigl(\sum_{i=1}^NT_i\bigr)\le 3\,L(c)$.
Hence, using Lemma \ref{bhu8y} again,
and the assumptions $L(c)< 1$ and $c< 1$,
\begin{align*}
d_W\Bigl(\frac{S}{\sigma}\,,\,Z\Bigr)&\le d_W\Bigl(\frac{S}{\sigma}\,,\,\sum_{i=1}^NV_i\Bigr)+d_W\bigl(\,\sum_{i=1}^NV_i,\,\sigma_c Z\bigr)+d_W(\sigma_c Z,\,Z)
\\
&\le E\Abs{\frac{S}{\sigma}-\sum_{i=1}^NV_i}+\sigma_c\,d_W\Bigl(\frac{\sum_{i=1}^NV_i}{\sigma_c}\,,\,Z\Bigr)+\abs{\sigma_c -1}
\\
&\le\sqrt{\text{Var}\,\bigl(\,\sum_{i=1}^NT_i\bigr)}+16\,\bigl(6\,c\,\sigma_c^2\bigr)^{1/3}+13\,L(c)
\\&
\le\sqrt{3\,L(c)}+16\,(6\,c)^{1/3}\,\bigl(1+13\,L(c)\bigr)^{2/3}+13\,L(c)
\\&\le\bigl(\sqrt{3}+13\bigr)\,L(c)^{1/2}+16\,(6\,c)^{1/3}\,\Bigl(1+\bigl(13\,L(c)\bigr)^{2/3}\Bigr)
\\&\le 16\,(6\,c)^{1/3}+\Bigl(\sqrt{3}+13+16\cdot 6^{1/3}\cdot (13)^{2/3}\Bigr)\,L(c)^{1/2}
\\&\le 30\,c^{1/3}+180\,L(c)^{1/2}.
\end{align*}
This concludes the proof of Theorem \ref{m1}.
\end{proof}

\medskip

\section{Total variation distance and rate of convergence}\label{tv21d}

Theorems \ref{math187v4} and \ref{m2} immediately imply the following result.

\begin{thm}\label{void56g7}
Let $\phi_n$ be the characteristic function of $S_n/\sigma_n$ and
$$l_n=2\,\int_0^\infty t\,\abs{\phi_n(t)}\,dt.$$
If conditions \eqref{c1}--\eqref{c2} hold, then
$$d_{TV}\Bigl(\frac{S_n}{\sigma_n},\,Z\Bigr)\le\sqrt{120}\,\Bigl\{c^{1/3}\kern-0.6mm+\kern-0.6mm12\,U_n(c/2)^{1/2}\Bigr\}^{1/2}
+\,30^{1/3}\,l_n^{2/3}\,\Bigl\{c^{1/3}+12\,U_n(c/2)^{1/2}\Bigr\}^{1/3}$$
for all $n\ge 1$ and $c>0$, where $Z$ is a standard normal random variable.
\end{thm}

\smallskip\noindent
\begin{proof}
First apply Theorem \ref{math187v4}, with $V\kern-0.3mm=\kern-0.3mm1$ and $X_n\kern-0.5mm=\kern-0.5mm \frac{S_n}{\sigma_n}$, and then use Theorem\,\ref{m2}.
\end{proof}

Obviously, Theorem \ref{void56g7} is non-trivial only if $l_n<\infty$. In this case, the probability distribution of $S_n$ is absolutely continuous. An useful special case is when conditions \eqref{c1}--\eqref{c2} hold and
\begin{gather}\label{m92dz6hjf55t}
\max_i\,\abs{X_{n,i}}\le\sigma_n\gamma_n\quad\quad\text{a.s. for some constants }\gamma_n.
\end{gather}
Under \eqref{m92dz6hjf55t}, since $U_n(m_n\,\gamma_n)=0$, Theorem
\ref{void56g7} yields
\begin{gather*}
d_{TV}\Bigl(\frac{S_n}{\sigma_n},\,Z\Bigr)\le \sqrt{120}\,\,(2\,m_n\,\gamma_n)^{1/6}+\,30^{1/3}\,l_n^{2/3}\,(2\,m_n\,\gamma_n)^{1/9}.
\end{gather*}
Sometimes, this inequality allows to obtain a CLT in total variation distance; see Example \ref{z5m8jjn} below.

\medskip

Finally, we discuss the convergence rate provided by Theorem \ref{m2} and we compare it with some existing results.

\medskip

A first remark is that Theorem \ref{m2} is calibrated to the dependence
case, and that it is not optimal in the independence case. To see this, it
suffices to recall that we assume $m_n\ge 1$ for all $n$. If
$X_{n1},\ldots,X_{nN_n}$ are independent, the best one can do is to let
$m_n=1$, but this choice of $m_n$ is not efficient as is shown by the following
example.

\begin{ex}
Suppose $X_{n1},\ldots,X_{nN_n}$ are independent and conditions \eqref{c2}
and \eqref{m92dz6hjf55t} hold. Define $m_n=1$ for all $n$. Then,
$U_n(\gamma_n)=0$ and Theorem \ref{m2} implies
$d_W\bigl(S_n/\sigma_n,\,Z\bigr)\le 30\,(2\,\gamma_n)^{1/3}$. However, the
Bikelis nonuniform inequality yields
$$\Abs{P(S_n/\sigma_n\le t)-P(Z\le
  t)}\le\frac{b}{(1+\abs{t})^3}\,\sum_{i=1}^{N_n}E\Bigl\{\frac{\abs{X_{n,i}}^3}{\sigma_n^3}\Bigr\}
\le\frac{b\,\gamma_n}{(1+\abs{t})^3}$$
for all $t\in\mathbb{R}$ and some universal constant $b$; see
e.g.\ \cite[p.~659]{DG}.
Hence,
$$d_W\Bigl(\frac{S_n}{\sigma_n},\,Z\Bigr)=\int_{-\infty}^\infty\abs{P(S_n/\sigma_n\le t)-P(Z\le t)}\,dt\le\int_{-\infty}^\infty\frac{b\,\gamma_n}{(1+\abs{t})^3}\,dt=b\,\gamma_n.$$
\end{ex}
\emph{}
\medskip

Leaving independence aside, a recent result to be mentioned is
\cite[Corollary 4.3]{DMR} by Dedecker, Merlevede and Rio. This result applies to sequences of random variables and requires a certain mixing condition (denoted by $(H_1)$) which is automatically true when $m_n=m$ for all $n$. In this case, under conditions \eqref{c2} and \eqref{m92dz6hjf55t}, one
obtains
\begin{gather}\label{ine22q}
d_W\Bigl(\frac{S_n}{\sigma_n},\,Z\Bigr)\le b\,\gamma_n\,\Bigl(1+c_n\,\log\,\bigl(1+c_n\,\sigma_n^2\bigr)\Bigr)
\end{gather}
where $b$ and $c_n$ are suitable constants with $b$ independent of $n$. Among other conditions, the $c_n$ must satisfy
$$c_n\,\sigma_n^2\ge\sum_{i=1}^{N_n}EX_{n,i}^2.$$
Inequality \eqref{ine22q} is actually sharp. However, if compared with
Theorem \ref{m2}, it has three drawbacks. First, unlike Theorem \ref{m2}, it
requires condition \eqref{m92dz6hjf55t}. Secondly, the mixing condition
$(H_1)$ is not easily verified unless $m_n=m$ for all $n$. Thirdly,
as seen in the next example,
even if \eqref{m92dz6hjf55t} holds and $m_n=m$ for all $n$, it may be that
$$\gamma_n\rightarrow 0\quad\text{but}\quad\gamma_n\,c_n\,\log\,\bigl(1+c_n\,\sigma_n^2\bigr)\rightarrow\infty\quad\quad\text{as }n\rightarrow\infty.$$
In such situations, Theorem \ref{m2} works while inequality \eqref{ine22q} does not.

\begin{ex}\label{tyt778e4}
Let $(a_n)$ be a sequence of numbers in $(0,1)$ such that $\lim_na_n=0$. Let $(T_i:i\ge 0)$ and $(V_{n,i}:n\ge 1,\,1\le i\le n)$ be two independent collections of real random variables. Suppose $(T_i)$ is i.i.d.\ with $P(T_0=\pm 1)=1/2$ and $V_{n,1},\ldots,V_{n,n}$ are i.i.d.\ with $V_{n,1}$ uniformly distributed on the set $(-1,-1+a_n)\cup (1-a_n,1)$.

\medskip

Fix a constant $\alpha\in (0,\,1/3)$ and define $N_n=n$ and
$$X_{n,i}=n^{-1/2}V_{n,i}+n^{-\alpha}(T_i-T_{i-1})$$
for $i=1,\ldots,n$. The array $(X_{n,i})$ is centered and 1-dependent (namely, $m_n=1$ for all $n$). In addition,
$S_n=n^{-1/2}\sum_{i=1}^nV_{n,i}+n^{-\alpha}(T_n-T_0)$ and $$\sigma_n^2=EV_{n,1}^2+2\,n^{-2\alpha},\quad\sum_{i=1}^nEX_{n,i}^2=EV_{n,1}^2+2\,n^{1-2\alpha}.$$
Since $\lim_n\sigma_n^2=\lim_n EV_{n,1}^2=1$, one obtains
$$\max_i\frac{\abs{X_{n,i}}}{\sigma_n}\le\frac{n^{-1/2}+2\,n^{-\alpha}}{\sigma_n}\le\frac{3\,n^{-\alpha}}{\sigma_n}<4\,n^{-\alpha}\quad\text{for large }n.$$
Hence, for large $n$, condition \eqref{m92dz6hjf55t} holds with
$\gamma_n=4\,n^{-\alpha}$. Since $U_n(4\,n^{-\alpha})=0$, Theorem \ref{m2} implies
(taking $c=8 n^{-\alpha}$)
$$d_W\Bigl(\frac{S_n}{\sigma_n},\,Z\Bigr)\le 60\,n^{-\alpha/3}\quad\quad\text{for large }n.$$
However,
\begin{align*}
4\,n^{-\alpha}\,c_n\,\log\,\bigl(1+c_n\,\sigma_n^2\bigr)&\ge 4\,n^{-\alpha}\,\frac{1}{\sigma_n^2}\,\sum_{i=1}^nEX_{n,i}^2\,\,\log\,\Bigl(1+\sum_{i=1}^nEX_{n,i}^2\Bigr)
\\&\ge 4\,(1-2\alpha)\,\,\frac{n^{1-3\alpha}}{\sigma_n^2}\,\,\log n\longrightarrow\infty.
\end{align*}
\end{ex}

\medskip

In addition to \cite[Corollary 4.3]{DMR}, there are some other estimates of $d_W\bigl(S_n/\sigma_n,\,Z\bigr)$. Without any claim of exhaustivity, we mention Fan and Ma \cite{FANMA}, R\"ollin \cite{ROLLIN} and Van Dung, Son and Tien \cite{VST} (R\"ollin's result has been used for proving Lemma \ref{vb1a0i}). There are also a number of estimates of $d_K\bigl(S_n/\sigma_n,\,Z\bigr)$ which, through Lemma \ref{b66z3e0ju5}, can be turned into upper bounds for $d_W\bigl(S_n/\sigma_n,\,Z\bigr)$; see \cite{DMR}, \cite{FANMA} and references therein. However, to our knowledge, none of these estimates implies Theorem \ref{m2}. Typically, they require further conditions (in addition to \eqref{c1}--\eqref{c2}) and/or they yield a worse convergence rate; see e.g.\ Remark \ref{hb88j9w} and Example \ref{tyt778e4}. This is the current state of the art. Our conjecture is that, under conditions \eqref{c1}--\eqref{c2} and possibly \eqref{m92dz6hjf55t}, the rate of Theorem \ref{m2} can be improved. To this end, one possibility could be using an upper bound provided by Haeusler and Joos \cite{HAEUJOOS} in the martingale CLT. Whether the rate of Theorem \ref{m2} can be improved, however, is currently an {\em open problem}.

\medskip

We conclude the paper with a CLT in total variation distance obtained via Theorem \ref{void56g7}.

\begin{ex}\label{z5m8jjn}
Let $(X_{n,i})$ and $(V_{n,i})$
be as in Example \ref{tyt778e4}. Denote by $\psi_n$ the characteristic
function of $\sum_{i=1}^nV_{n,i}$.
Then, for each $t\in\mathbb{R}$,
\begin{gather*}
\psi_n(t)=\Bigl(\,\frac{1}{a_n}\,\int_{1-a_n}^1\cos(t\,x)\,dx\Bigr)^n\quad\quad\text{ and}
\\\abs{\phi_n(t)}\le\Abs{\psi_n\bigl[t\,(n\,\sigma_n^2)^{-1/2}\bigr]}=\Abs{\,\frac{1}{a_n}\,\int_{1-a_n}^1\cos\bigl[t\,(n\,\sigma_n^2)^{-1/2}\,x\bigr]\,dx}^n.
\end{gather*}
After some algebra (we omit the explicit calculations) it can be shown that
$$l_n=2\,\int_0^\infty t\,\abs{\phi_n(t)}\,dt\le b\,a_n^{-2}$$
for some constant $b$ independent of $n$. Recalling that $m_n=1$ and
$\gamma_n=4\,n^{-\alpha}$ for large $n$ (see Example \ref{tyt778e4}),
Theorem \ref{void56g7} yields (taking again $c=2m_n\gamma_n=8 n^{-\alpha}$)
\begin{align*}
d_{TV}\Bigl(\frac{S_n}{\sigma_n},\,Z\Bigr)&\le\sqrt{120}\,\,(2\,m_n\,\gamma_n)^{1/6}+\,30^{1/3}\,l_n^{2/3}\,(2\,m_n\,\gamma_n)^{1/9}
\\&\le\sqrt{120}\,\,8^{1/6}\,n^{-\alpha/6}\,+\,30^{1/3}\,b^{2/3}\,8^{1/9}\,\Bigl(a_n^4\,n^{\alpha/3}\bigr)^{-1/3}
\end{align*}
for large $n$. Therefore, the probability distribution of $S_n/\sigma_n$ converges to the standard normal law, in total variation distance, provided $a_n^4\,n^{\alpha/3}\rightarrow\infty$.
\end{ex}


\begin{thebibliography}{9}

\bibitem{BERG} Bergstr\"om H. (1970) A comparison method for distribution
  functions of sums of independent and dependent random variables,
  (Russian), {\em Teor. Verojatnost. i Primenen} 15, 442--468, 750;
English transl. {\em Theor. Probab. Appl.} 15 (1970), 430--457, 727.
MR 0283850, MR 0281245

\bibitem{BERK} Berk K.N. (1973) A central limit theorem for $m$-dependent random variables with unbounded $m$, {\em Ann. Probab.} 1, 352--354. MR 0350815

\bibitem{BRAD}  Bradley R.C. (2007) {\em Introduction to strong mixing conditions}, Vol. 1--3, Kendrick Press, Heber City, UT. MR 2325294--2325296


\bibitem{CHEN}
Chen L.H.Y.,  Goldstein L.,  Shao Q.-M. (2011)
\emph{Normal approximation by Stein's method},
Springer-Verlag, Berlin.
MR 2732624

\bibitem{DG} DasGupta A. (2008) {\em Asymptotic theory of statistics and probability}, Springer, New York.

\bibitem{DMR} Dedecker J., Merlevede F., Rio E. (2021) Rates of convergence in the central limit theorem for martingales in the non stationary setting, {\em arXiv:2101.06956v1}.

\bibitem{DIAN} Diananda P.H. (1955) The central limit theorem for $m$-dependent variables, {\em Proc. Cambridge Philos. Soc.} 51, 92--95. MR 0067396

\bibitem{FANMA} Fan X., Ma X. (2020) On the Wasserstein distance for a
  martingale central limit theorem, {\em Statist. Probab. Lett.} 167,  108892.
MR 4138415

\bibitem{HAEUJOOS} Haeusler E., Joos K. (1988) A nonuniform bound on the rate of convergence in the martingale central limit theorem, {\em Ann. Probab.} 16, 1699--1720.

\bibitem{HB} Heyde C.C., Brown B.M. (1970) On the departure from normality of a certain class of martingales, {\em Ann. Math. Stat.} 41, 2161--2165.
MR 0293702

\bibitem{HOEROB} Hoeffding W., Robbins H. (1948) The central limit theorem for dependent random variables, {\em Duke Math. J.} 15, 773--780. MR 0026771


\bibitem{OREY} Orey S. (1958) A central limit theorem for $m$-dependent random variables, {\em Duke Math. J.} 25, 543--546. MR 0097841

\bibitem{PEL} Peligrad M. (1996) On the asymptotic normality of sequences of weak dependent random variables, {\em J. Theoret. Probab.} 9, 703--715. MR 1400595

\bibitem{PR} Pratelli L., Rigo P. (2018) Convergence in total variation to a mixture of Gaussian laws, {\em Mathematics} 6, 99.

\bibitem{RIO} Rio E. (1995) About the Lindeberg method for strongly mixing sequences, {\em ESAIM: Probab. Statist.} 1, 35--61. MR 1382517

\bibitem{ROLLIN} R\"ollin A. (2018) On quantitative bounds in the mean
  martingale central limit theorem, {\em Statist. Probab. Lett.} 138, 171--176.
MR 3788734

\bibitem{RW} Romano J.P., Wolf M. (2000) A more general central limit theorem for $m$-dependent random variables with unbounded $m$, {\em Statist. Probab. Lett.} 47, 115--124. MR 1747098

\bibitem{UT1} Utev S.A. (1987)
On the central limit theorem for $\varphi$-mixing arrays of random variables,
{\em Theor. Probab. Appl.} 35, 131--139.
MR 1050059

\bibitem{UT2} Utev S.A. (1990) Central limit theorem for dependent random
  variables,
{\em Probability theory and mathematical statistics, Vol. II} (Vilnius, 1989),
519--528, Mokslas, Vilnius.
MR 1153906

\bibitem{VST} Van Dung L., Son T.C., Tien N.D. (2014) $L_1$-Bounds for some martingale central limit theorems, {\em Lith. Math. J.} 54, 48--60.
MR 3189136








\end{thebibliography}
\end{document}